\documentclass[10pt]{article}
\usepackage{mathrsfs}
\usepackage{amsthm}
\usepackage{amsmath}
\usepackage{graphicx}
\usepackage{color}
\usepackage{amsfonts}
\newtheorem{theorem}{Theorem}[section]

\newtheorem{lemma}[theorem]{Lemma}

\numberwithin{equation}{section}
\normalsize

\begin{document}
\title{\textbf{Convergence rates for subcritical threshold-one contact processes on lattices}}
\author{Xiaofeng Xue \thanks{\textbf{E-mail}: xuexiaofeng@ucas.ac.cn \textbf{Address}: School of Mathematical Sciences, University of Chinese Academy of Sciences, Beijing 100049, China.}\\ University of Chinese Academy of Sciences}

\date{}
\maketitle

\noindent {\bf Abstract:}

In this paper we are concerned with threshold-one contact processes
on lattices. We show that the probability that the origin is
infected converges to $0$ at an exponential rate $I$ in the
subcritical case. Furthermore, we give a limit theorem for $I$ as
the degree of the lattice grows to infinity. Our results also hold
for classic contact processes on lattices.


\section{Introduction}\label{section 1 introduction}
In this paper we are concerned with threshold-one contact processes on lattices $\mathbb{Z}^d,d=1,2,\ldots$
For any $x,y\in \mathbb{Z}^d$, we write $x\sim y$ when there is an edge connecting these two vertices.
We say that $x$ and $y$ are neighbors when $x\sim y$.

The threshold-one contact process $\{\eta_t\}_{t\geq 0}$ on $\mathbb{Z}^d$ is with state space $\{0,1\}^{\mathbb{Z}^d}$.
In other words, at each vertex of $\mathbb{Z}^d$ there is a spin taking value $0$ or $1$.
For each $x\in \mathbb{Z}^d$ and $t>0$, the spin at $x$ at moment $t$ is denoted by $\eta_t(x)$. Furthermore, we define $\eta_{t-}(x)$ as
\[
\eta_{t-}(x):=\lim_{s<t,s\uparrow t}\eta_s(x).
\]
Hence $\eta_{t-}(x)$ is the spin at $x$ at the moment just before $t$.

$\{\eta_t\}_{t\geq 0}$ evolves according to independent Poisson
processes $\{N_x(t):t\geq 0\}_{x\in \mathbb{Z}^d}$ and
$\{Y_x(t):t\geq 0\}_{x\in \mathbb{Z}^d}$. For each $x\in
\mathbb{Z}^d$, $N_x$ is with rate $1$ and $Y_x$ is with rate
$\lambda$, where $\lambda>0$ is a parameter called the infection
rate. At $t=0$, each spin takes a value from $\{0,1\}$ according to
some probability distribution. Then, for each $x\in \mathbb{Z}^d$,
the spin at $x$ may flip only at event times of $N_x$ and $Y_x$. For
any event time $s$ of $N_x$, $\eta_s(x)=0$ no matter whatever
$\eta_{s-}(x)$ is. For any event time $r$ of $Y_x$, if
$\eta_{r-}(x)=1$, then $\eta_{r}(x)=1$. If $\eta_{r-}(x)=0$, then
$\eta_{r}(x)=1$ when and only when at least one neighbor $y$ of $x$
satisfies $\eta_{r-}(y)=1$.

Therefore, $\{\eta_t\}_{t\geq 0}$ is a spin system (see Chapter 3 of \cite{LIG1985}) with flip rates function given by
\begin{equation}\label{equ 1.1 rate function}
c(x,\eta)=
\begin{cases}
1&\text{~if~}\eta(x)=1,\\
\lambda &\text{~if~} \eta(x)=0\text{~and~}\sum_{y:y\sim x}\eta(y)\geq 1,\\
0&\text{otherwise}
\end{cases}
\end{equation}
for any $(x,\eta)\in \mathbb{Z}^d\times \{0,1\}^{\mathbb{Z}^d}$.

Intuitively, the threshold-one contact process describes the spread of an infected disease. Vertices with spin $1$ are infected individuals while vertices with spin $0$ are healthy. An infected vertex waits for an exponential time with rate one to become healthy while a healthy vertex is infected by neighbors with rate $\lambda$ when at least one neighbor is infected.

Our main result in this paper about the threshold-one contact process $\{\eta_t\}_{t\geq 0}$ also holds for the classic contact process $\{\beta_t\}_{t\geq 0}$. The flip rates function of $\beta_t$ is given by
\begin{equation}\label{equ 1.2}
\widehat{c}(x,\beta)=
\begin{cases}
1&\text{~if~} \beta(x)=1,\\
\lambda\sum_{y:y\sim x}\eta(y)&\text{~if~} \beta(x)=0
\end{cases}
\end{equation}
for any $(x,\beta)\in \mathbb{Z}^d\times \{0,1\}^{\mathbb{Z}^d}$. The main difference between $\eta_t$ and $\beta_t$ is that for $\beta_t$, a healthy vertex is infected at rate proportional to the number of infected neighbors.

The threshold-one contact processes is introduced in \cite{Dur1991} by Cox and Durrett as a tool to study threshold voter models (see Part two of \cite{LIG1999} and \cite{And1992, Handjani1999, LIG1994, Xue2012, Xue2015}). \cite{Dur1991} gives an important dual relationship between the threshold-one contact process and an additive Markov processes. According to this dual relationship, \cite{Dur1991} shows that the critical value $\lambda_c(d)$ for the threshold-one contact process on $\mathbb{Z}^d$ satisfies $\lambda_c(d)\leq 2.18/d$. \cite{Xue2014} develops this result by showing that $\lim_{d\rightarrow+\infty}2d\lambda_c(d)=1$.
In recent years, there are some works concerned with threshold contact processes with threshold $K>1$. \cite{Mou2009} shows that the critical value   $\lambda_c(d,K)$ for the threshold $K>1$ contact process on $\mathbb{Z}^d$ satisfies $\lim_{d\rightarrow+\infty}\lambda_c(d,K)=0$. \cite{Fonte2008} shows that the same conclusion holds for the case on regular trees $\mathbb{T}^N$ and gives the rate at which $\lambda_c(\mathbb{T}^N,K)$ converges to $0$ as $N$ grows to infinity.

\section{Main result}\label{section 2 main result}
In this section, we will give the main result of this paper. First
we introduce some notations. For $d\geq 1$ and $\lambda>0$, we
denote by $P_{\lambda,d}$ the probability measure of the
threshold-one contact process $\{\eta_t\}_{t\geq 0}$ on
$\mathbb{Z}^d$ with infection rate $\lambda$. We denote by
$E_{\lambda,d}$ the expectation operator with respect to
$P_{\lambda,d}$. We write $\eta_t$ as $\eta_t^{\eta}$ when
\[P_{\lambda,d}(\eta_0=\eta)=1\]
for some $\eta\in \{0,1\}^{\mathbb{Z}^d}$. We denote by $\delta_1$ and $\delta_0$ configurations in $\{0,1\}^{\mathbb{Z}^d}$ such that
\[
\delta_1(x)=1,~\delta_0(x)=0
\]
for each $x\in \mathbb{Z}^d$. We denote by $O$ the origin of $\mathbb{Z}^d$ and denote by $e_1$ the unit vector $(1,0,0,\ldots,0)$.

Since the threshold-one contact process is attractive (see Chapter 3 of \cite{LIG1985}), for any $t>s$ and $\lambda_1>\lambda_2$,
\[
P_{\lambda_1,d}(\eta_s^{\delta_1}(O)=1)\geq P_{\lambda_2,d}(\eta_t^{\delta_1}(O)=1).
\]
As a result, it is reasonable to define the following critical value.
\begin{equation}\label{equ 2.3}
\lambda_c(d):=\sup\{\lambda:\lim_{t\rightarrow+\infty}P_{\lambda,d}(\eta_t^{\delta_1}(O)=1)=0\}
\end{equation}
for $d\geq 1$.

When $\lambda<\lambda_c(d)$, the process $\eta_t$ converges weakly to $\delta_0$ as $t\rightarrow+\infty$, which is called the subcritical case.

In the subcritical case, we are concerned with the rate at which the
probability that $O$ is infected converges to $0$ as the time $t$
grows to infinity. To introduce our main result, we give a lemma at
first.
\begin{lemma}\label{lemma 2.1}
For any $\lambda\geq0$ and $d\geq 1$, there exists $I(\lambda,d)\in [-\infty,0]$ such that
\begin{equation}\label{equ 2.1}
\lim_{t\rightarrow+\infty}\frac{1}{t}\log P_{\lambda,d}(\eta^{\delta_1}_t(O)=1)=I(\lambda,d).
\end{equation}
\end{lemma}
After giving $\lambda$ a proper scale, we obtain the following limit theorem of $I(\lambda,d)$ as our main result.
\begin{theorem}\label{theorem 2.2 main}
For any $\lambda\geq 0$,
\begin{equation}\label{equ 2.2}
\lim_{d\rightarrow+\infty}I\big(\frac{\lambda}{d},~d\big)=
\begin{cases}
2\lambda-1 & \text{~if~} \lambda\in [0,1/2],\\
0 &\text{~if~}\lambda>1/2.
\end{cases}
\end{equation}
\end{theorem}
Theorem \ref{theorem 2.2 main} shows that for  subcritical
threshold-one contact process with infection rate $\lambda$, the
probability that $O$ is infected converges to $0$ as $t\rightarrow
+\infty$ at an exponential rate approximate to $2\lambda d-1$ when
the dimension $d$ is sufficiently large.

According to the dual relationship introduced in \cite{Dur1991},
there is an intuitive explanation for Theorem \ref{theorem 2.2
main}. When the dimension $d$ is large, the threshold-one contact
process is similar with a branching process such that each particle
generates $2d$ particles at rate $\lambda$ or dies at rate one. The
mean of the sum of the particles at $t$ is $\exp\{(2\lambda
d-1)t\}$.

In Theorem \ref{theorem 2.2 main}, the case where $\lambda>1/2$ is trivial, since \cite{Xue2014} shows that
\[
\lim_{d\rightarrow+\infty}2d\lambda_c(d)=1.
\]

Similar conclusion with Theorem \ref{theorem 2.2 main} holds for the classic contact process $\{\beta_t\}_{t\geq 0}$, the flip rate function of which is given in \eqref{equ 1.2}.

\begin{theorem}\label{theorem 2.3}
For any $\lambda\geq 0$ and $d\geq 1$, there exists $J(\lambda,d)\in[-\infty,0]$ such that
\[
\lim_{t\rightarrow+\infty}\frac{1}{t}\log P_{\lambda,d}(\beta_t^{\delta_1}(O)=1)=J(\lambda,d)
\]
and
\[
\lim_{d\rightarrow+\infty}J(\frac{\lambda}{d},~d)=
\begin{cases}
2\lambda-1 &\text{~if~} \lambda\in [0,1/2],\\
0 &\text{~if~} \lambda>1/2.
\end{cases}
\]
\end{theorem}

In this paper, the proof of theorem about $\beta_t$ is similar with that of the counterpart conclusion about $\eta_t$. We will give all the details in the proof of theorem about $\eta_t$ and give just a sketch for the proof of theorem about $\beta_t$.

At the end of this section, we give the proof of Lemma \ref{lemma 2.1}. The proof of Theorem \ref{theorem 2.2 main} is divide into Section \ref{section 3 upper bound} and Section \ref{section 4 lower bound}.

\proof[Proof of Lemma \ref{lemma 2.1}]
We utilize the dual relationship introduced in \cite{Dur1991}.
Let $\{A_t\}_{t\geq 0}$ be a Markov process with state space
\[
2^{\mathbb{Z}^d}:=\{A:A\subseteq \mathbb{Z}^d\}
\]
and flip rate functions
\[
A_t\rightarrow
\begin{cases}
A_t\setminus x & \text{~at rate~} 1,\\
A_t\cup \{y:y\sim x\}& \text{~at rate~} \lambda \\
\end{cases}
\]
for any $t\geq 0$ and each $x\in A_t$.

We write $A_t$ as $A_t^A$ when $A_0=A$. Then, according to \cite{Dur1991}, there is a dual relationship between $\{\eta_t\}_{t\geq 0}$ and $\{\beta_t\}_{t\geq 0}$ such that
\begin{equation}\label{equ 2.4}
P_{\lambda,d}(\eta_t^{\delta_1}(O)=1)=P_{\lambda,d}(A_t^O\neq \emptyset).
\end{equation}

As a result, according to strong Markov property,
\begin{align}\label{equ 2.5}
P_{\lambda,d}(\eta_{t+s}^{\delta_1}(O)=1)=P_{\lambda,d}(A_{t+s}^O\neq \emptyset)=E_{\lambda,d}\Big[P_{\lambda,d}\big(A_s^{A_t}\neq \emptyset\big);A_t^O\neq \emptyset\Big].
\end{align}
Since $A_t$ is symmetric for $Z^d$ and is a monotone process under the partial order that $A\geq B$ if and only if $A\supseteq B$,
\begin{equation}\label{equ 2.6}
P_{\lambda,d}\big(A_s^{A_t}\neq\emptyset \big)\geq P_{\lambda,d}(A_s^O\neq \emptyset)
\end{equation}
on the event $\{A_t\neq \emptyset\}$.

By \eqref{equ 2.5} and \eqref{equ 2.6},
\begin{align*}
P_{\lambda,d}(\eta_{t+s}^{\delta_1}(O)=1)&\geq P_{\lambda,d}(A_t^O\neq \emptyset)P_{\lambda,d}(A_s^O\neq \emptyset)\\
&=P_{\lambda,d}(\eta_{t}^{\delta_1}(O)=1)P_{\lambda,d}(\eta_{s}^{\delta_1}(O)=1)
\end{align*}
and hence,
\begin{equation}\label{equ 2.7}
\log P_{\lambda,d}(\eta_{t+s}^{\delta_1}(O)=1)\geq \log P_{\lambda,d}(\eta_{t}^{\delta_1}(O)=1)+\log P_{\lambda,d}(\eta_{s}^{\delta_1}(O)=1)
\end{equation}
for any $t,s\geq 0$.

The existence of $I(\lambda,d)$ follows from \eqref{equ 2.7} and Fekete's Subadditive Lemma. By Fekete's Subadditive Lemma,
\[
I(\lambda,d)=\sup_{t\geq 0}\frac{1}{t}\log P_{\lambda,d}(\eta_t^{\delta_1}(O)=1).
\]

\qed

The proof of the existence of $J(\lambda,d)$ is nearly the same as that of $I(\lambda,d)$ by the self-duality of $\{\beta_t\}_{t\geq 0}$ introduced in Theorem 6.1.7 of \cite{LIG1985}.

\proof[Proof of the existence of $J(\lambda,d)$]

Let $C_t=\{x\in Z^d:\beta_t(x)=1\}$ and write $C_t$ as $C_t^A$ when $C_0=A$, then according to Theorem 6.1.7 of \cite{LIG1985},
\begin{equation}\label{equ 2.8}
P_{\lambda,d}(\beta^{\delta_1}_t(O)=1)=P_{\lambda,d}(C_t^O\neq \emptyset).
\end{equation}
The existence of $J(\lambda,d)$ follows from \eqref{equ 2.8} and a
similar analysis with that after \eqref{equ 2.4} in the proof of
Lemma \ref{lemma 2.1}.

\qed

\section{Upper bound}\label{section 3 upper bound}
In this section we will give upper bounds for $I(\lambda/d,~d)$ and $J(\lambda/d,~d)$.

The proofs of Theorem \ref{theorem 2.2 main} for cases where $\lambda=0$ and $\lambda>1/2$ are trivial. According to \cite{Xue2014},
\[
\lim_{d\rightarrow+\infty}2d\lambda_c(d)=1.
\]
As a result, for $\lambda>1/2$ and sufficiently large $d$, $\lambda_c(d)<\lambda/d$ and
\[
\lim_{t\rightarrow+\infty}P_{\lambda/d,d}(\eta_t^{\delta_1}(O)=1)=K(\lambda,d)>0.
\]
Therefore,
\[
I(\lambda/d,d)=\lim_{t\rightarrow+\infty}\frac{1}{t}\log P_{\lambda/d,d}(\eta_t^{\delta_1}(O)=1)=0
\]
for $\lambda>1/2$ and sufficiently large $d$.

The above analysis also holds for $J(\lambda/d,~d)$ since \cite{Grif1983} shows that the critical value $\widehat{\lambda}_c(d)$ for the classic contact process $\{\beta_t\}_{t\geq 0}$ on $Z^d$ satisfies
\[
\lim_{d\rightarrow+\infty}2d\widehat{\lambda}_c(d)=1.
\]

When $\lambda=0$, $O$ waits for an exponential time with rate one to become healthy and will never be infected again. Hence,
\[
P_{0,d}(\eta_t^{\delta_1}(O)=1)=P_{0,d}(\beta_t^{\delta_1}(O)=1)=e^{-t}
\]
and
\[
I(0,d)=J(0,d)=\lim_{t\rightarrow+\infty}\log e^{-t}=-1.
\]

Now we only need to deal with the case where $\lambda\in (0,1/2)$.
The following lemma gives an upper bound for $I(\lambda,d)$.
\begin{lemma}\label{lemma 3.1}
For any $\lambda>0$ and $d\geq 1$,
\[
\max\{I(\lambda,d),~J(\lambda,d)\}\leq 2\lambda d-1.
\]
\end{lemma}
As a direct corollary,
\[
\max\big\{\limsup_{d\rightarrow+\infty}I(\lambda/d,d),~\limsup_{d\rightarrow+\infty}J(\lambda/d,d)\big\}\leq 2\lambda-1
\]
for $\lambda\in (0,1/2)$.

\proof[Proof of Lemma \ref{lemma 3.1}]

According to the flip rate functions of $\{\eta_t\}_{t\geq 0}$ given in \eqref{equ 1.1 rate function} and Hille-Yosida Theorem,
\begin{align*}
\frac{d}{dt}P_{\lambda,d}(\eta_t^{\delta_1}(O)=1)=&-P_{\lambda,d}(\eta_t^{\delta_1}(O)=1)\\
&+\lambda P_{\lambda,d}(\eta_t^{\delta_1}(O)=0,\exists~y\sim O, \eta_t^{\delta_1}(y)=1)\\
&\leq -P_{\lambda,d}(\eta_t^{\delta_1}(O)=1)+\lambda P_{\lambda,d}(\exists~y\sim O, \eta_t^{\delta_1}(y)=1)\\
&\leq -P_{\lambda,d}(\eta_t^{\delta_1}(O)=1)+\lambda \sum_{y:y\sim O}P_{\lambda,d}(\eta_t^{\delta_1}(y)=1)\\
&=(2\lambda d-1)P_{\lambda,d}(\eta_t^{\delta_1}(O)=1),
\end{align*}
since $O$ has $2d$ neighbors and $\{\eta_t\}_{t\geq 0}$ is symmetric for $Z^d$.

Then, according to Gr\"{o}nwall's inequality,
\[
P_{\lambda,d}(\eta_t^{\delta_1}(O)=1)\leq e^{(2\lambda d-1)t}P_{\lambda,d}(\eta_0^{\delta_1}(O)=1)=e^{(2\lambda d-1)t}
\]
and hence
\[
I(\lambda,d)=\lim_{t\rightarrow+\infty}\frac{1}{t}\log P_{\lambda,d}(\eta_t^{\delta_1}(O)=1)\leq 2\lambda d-1.
\]
The analysis for $J(\lambda,d)$ is similar. According to the flip rate functions given in \eqref{equ 1.2},
\begin{align*}
\frac{d}{dt}P_{\lambda,d}(\beta_t^{\delta_1}(O)=1)&=-P_{\lambda,d}(\beta_t^{\delta_1}(O)=1)\\
&+\lambda\sum_{y:y\sim O}P_{\lambda,d}(\beta_t^{\delta_1}(O)=0,\beta_t^{\delta_1}(y)=1)\\
&\leq -P_{\lambda,d}(\beta_t^{\delta_1}(O)=1)+\lambda\sum_{y:y\sim O}P_{\lambda,d}(\beta_t^{\delta_1}(y)=1).
\end{align*}
Then $J(\lambda,d)\leq 2\lambda d-1$ follows from the same analysis as that of $I(\lambda,d)$.

\qed

\section{Lower bound}\label{section 4 lower bound}
In this section we will give lower bounds for $I(\lambda/d,~d)$ and $J(\lambda/d,~d)$ for $\lambda\in (0,1/2)$. The main tool we utilize is a Markov process $\{\zeta_t\}_{t\geq 0}$ with state space $[0,+\infty)^{\mathbb{Z}^d}$ introduced in \cite{Xue2014}. In other words, for $\{\zeta_t\}_{t\geq 0}$, there is a spin at each vertex of $\mathbb{Z}^d$ taking a nonnegative value.

Let $\{N_x(t):t\geq 0\}_{x\in \mathbb{Z}^d}$ and $\{Y_x(t):t\geq
0\}_{x\in \mathbb{Z}^d}$ be the same Poisson processes as that in
Section \ref{section 1 introduction}. $\{\zeta_t\}_{t\geq 0}$
evolves according to $\{N_x\}_{x\in \mathbb{Z}^d}$ and
$\{Y_x\}_{x\in \mathbb{Z}^d}$. At $t=0$, $\zeta_0(x)>0$ for each
$x\in \mathbb{Z}^d$. For any event time $s$ of $N_x$, $\zeta_s(x)=0$
no matter whatever $\zeta_{s-}(x)$ is. For any event time $r$ of
$Y_x$, $\zeta_r(x)=\zeta_{r-}(x)+\sum_{y:y\sim x}\zeta_{r-}(y)$.
Between any adjacent event times of Poisson processes, $\zeta_t(x)$
evolves according to deterministic ODE
\[
\frac{d}{dt}\zeta_t(x)=(1-2\lambda d)\zeta_t(x).
\]
In other words, if there is no event time of $N_x$ or $Y_x$ in $[t_1,t_2]$, then
\begin{equation}\label{equ ODE 4.1}
\zeta_{t_2}(x)=\zeta_{t_1}(x)\exp\{(1-2\lambda d)(t_2-t_1)\}.
\end{equation}

It is useful for us to give the generator of $\{\zeta_t\}_{t\geq 0}$. For any $\zeta\in [0,+\infty)^{\mathbb{Z}^d}$, $x\in \mathbb{Z}^d$ and
$m\in [0,+\infty)$, we define $U(\zeta,x)=\zeta(x)+\sum_{y:y\sim x}\zeta(y)$ and define $\zeta^{x,m}\in [0,+\infty)^{\mathbb{Z}^d}$ as
\[
\zeta^{x,m}(y)=
\begin{cases}
\zeta(y)&\text{~if~}y\neq x,\\
m &\text{~if~}y=x.
\end{cases}
\]

Then, the generator $\Omega$ of $\{\zeta_t\}_{t\geq 0}$ is given by
\begin{align}\label{equ 4.1 generator}
\Omega f(\zeta)&=\sum_{x\in \mathbb{Z}^d}\big[f(\zeta^{x,0})-f(\zeta)\big]+\lambda\sum_{x\in \mathbb{Z}^d}\big[f(\zeta^{x,U(\zeta,x)})-f(\zeta)\big]\notag\\
&+(1-2\lambda d)\sum_{x\in \mathbb{Z}^d}f_x^{\prime}(\zeta)\zeta(x)
\end{align}
for $f\in C([0,+\infty)^{\mathbb{Z}^d})$, where $f^{\prime}_x(\zeta)$ is the partial derivative of $f(\zeta)$ with respect to the coordinate $\zeta(x)$.

The following lemmas are crucial for us to give lower bound for
$I(\lambda,d)$.
\begin{lemma}\label{lemma 4.1}
There is a coupling of $\eta_t^{\delta_1}$ and $\zeta_t$ such that
\[
\eta_t^{\delta_1}(x)=
\begin{cases}
1 &\text{~if~}\zeta_t(x)>0,\\
0 &\text{~if~}\zeta_t(x)=0
\end{cases}
\]
for each $x\in \mathbb{Z}^d$ and $t\geq 0$.
\end{lemma}

\proof

For any $t\geq 0$ and $x\in \mathbb{Z}^d$, we define
\[
\widetilde{\eta}_t(x) = \begin{cases} 1&\text{~if~} \zeta_t(x)>0,\\
0&\text{~if~} \zeta_t(x)=0.
\end{cases}
\]
Then, $\widetilde{\eta}_0=\delta_1$. At any event time $s$ of $N_x$,
$\zeta_s(x)=0$ and hence $\widetilde{\eta}_s(x)=0$. At any event
time $r$ of $Y_x$, $\widetilde{\eta}(x)$ flips from $0$ to $1$ if
and only if $\zeta_{r}(x)=0+\sum_{y:y\sim x}\zeta_{r-}(y)>0$. In
other words, conditioned on $\widetilde{\eta}_{r-}(x)=0$,
$\widetilde{\eta}_{r}(x)=1$ if and only if at least one neighbor $y$
of $x$ satisfies $\zeta_{r-}(y)>0$ and meanwhile
$\widetilde{\eta}_{r-}(y)=1$. According to ODE \eqref{equ ODE 4.1},
between any adjacent event times of Poisson processes $N_x$ and
$Y_x$, $\zeta_t(x)$ can not flip from positive value to zero or flip
from zero to positive value, which makes $\widetilde{\eta}_t(x)$
still.

Therefore, $\{\widetilde{\eta}_t\}_{t\geq 0}$ evolves in the same
way as that of $\{\eta_t\}_{t\geq 0}$. Since
$\widetilde{\eta}_0=\eta_0^{\delta_1}=\delta_1$,
$\{\widetilde{\eta}_t\}_{t\geq 0}$ and $\{\eta_t^{\delta_1}\}_{t\geq
0}$ have the same distribution.

\qed

\begin{lemma}\label{lemma 4.2}
When $\zeta_0=\zeta$ where $\zeta(x)>0$ for each $x\in
\mathbb{Z}^d$, then there exists $C(\lambda,d, \zeta)>0$ such that
\[
E_{\lambda,d}\zeta_t(O)\geq C(\lambda,d,\zeta)t^{-\frac{d}{2}}
\]
for any $t\geq 0$.
\end{lemma}

\proof

According to the generator of $\{\zeta_t\}_{t\geq 0}$ given in
\eqref{equ 4.1 generator} and Theorem 9.1.27 of \cite{LIG1985},
\begin{align}\label{equ 4.2 three}
\frac{d}{dt}E_{\lambda,d}\zeta_t(x)&=-E_{\lambda,d}\zeta_t(x)+\lambda\sum_{y:y\sim
x}E_{\lambda,d}\zeta_t(y)+(1-2\lambda d)E_{\lambda,d}\zeta_t(x)\notag\\
&=\lambda\sum_{y:y\sim x}E_{\lambda,d}\zeta_t(y)-2\lambda d
E_{\lambda,d}\zeta_t(x)
\end{align}
for each $x\in \mathbb{Z}^d$.

Let $Q=\big(q(x,y)\big)_{x,y\in \mathbb{Z}^d}$ be the Q-matrix of
the continuous time simple random walk on $Z^d$ such that
\[
q(x,y)=
\begin{cases}
\lambda & \text{~if~} y\sim x,\\
-2\lambda d &\text{~if~} y=x,\\
0 &\text{~else},
\end{cases}
\]
then, by \eqref{equ 4.2 three},
\[
E_{\lambda,d}\zeta_t=P_t\zeta_0,
\]
where
\[
P_t=\big(p_t(x,y)\big)_{x,y\in
\mathbb{Z}^d}=e^{tQ}=\sum_{n=0}^{+\infty}\frac{t^nQ^n}{n!}.
\]
In other words, $P_t$ is the transition function of the simple
random walk with Q-matrix $Q$.

According to classic theory of continuous time simple random walk on
$\mathbb{Z}^d$, there exists $C>0$ such that
\[
p_t(O,O)\geq [C(\lambda t)^{-\frac{1}{2}}]^d
\]
for any $t\geq 0$, where $C$ does not depend on $\lambda$ and $d$.

Therefore,
\begin{align*}
E_{\lambda,d}\zeta_t(O)&=\sum_{x\in
\mathbb{Z}^d}p_t(O,x)\zeta_0(x)\geq p_t(O,O)\zeta_0(O)\notag\\
&\geq\zeta_0(O)[C(\lambda t)^{-\frac{1}{2}}]^d.
\end{align*}
Let $C(\lambda,d,\zeta)=\zeta_0(O)C^d\lambda^d$, then the proof is
completed.

\qed

We define $F:[0,+\infty)\rightarrow [0,+\infty)^{\mathbb{}Z^d}$ as
\begin{equation}\label{equ 4.2}
F_t(x)=E_{\lambda ,d}\big[\zeta_t(O)\zeta_t(x)\big]
\end{equation}
for any $t\geq 0$ and $x\in \mathbb{Z}^d$. Then, the following lemma holds for $F$.

\begin{lemma}\label{lemma 4.3}
For any $t\geq 0$,
\begin{equation}\label{equ 4.2 two}
\frac{d}{dt}F_t=GF_t,
\end{equation}
where $G$ is a $\mathbb{Z}^d\times \mathbb{Z}^d$ matrix such that
\begin{equation}\label{equ 4.3}
G(x,y)=
\begin{cases}
-4\lambda d &\text{~if~} x=y\text{~and~}x\neq O,\\
2\lambda &\text{~if~}  y\sim x \text{~and~}x\neq O,\\
1-2\lambda d &\text{~if~} x=y=O,\\
2\lambda d &\text{~if~} x=O\text{~and~}y=e_1,\\
2\lambda d &\text{~if~} x=O,~z\sim e_1 \text{~and~}z\neq O, \\
0 &\text{~otherwise,}
\end{cases}
\end{equation}
where $e_1=(1,0,0,\ldots,0)$.
\end{lemma}

\proof

\eqref{equ 4.2 two} follows directly from the generator of
$\{\zeta_t\}_{t\geq}$ given in \eqref{equ 4.1 generator}
 and Theorem 9.3.1 of \cite{LIG1985}.

\qed

\begin{lemma}\label{lemma 4.4}
If $\lambda$ satisfies
\[
GH=\mu H
\]
for some $\mu>0$ and $H: \mathbb{Z}^d\rightarrow R$ such that
$H(x)>0$ for each $x\in \mathbb{Z}^d$, then
\[
I(\lambda,d)\geq -\mu.
\]

\end{lemma}

\proof Let $\zeta_0(x)=H(x)$ for each $x\in \mathbb{Z}^d$, then,
according to Lemma \ref{lemma 4.1}, Lemma \ref{lemma 4.2}, and
H\"{o}lder's inequality,
\begin{align}\label{equ 4.6}
P_{\lambda,d}(\eta_{t}^{\delta_1}(O)=1)&=P_{\lambda,d}(\zeta_t(O)>0)\notag\\
&\geq
\frac{[E_{\lambda,d}\zeta_t(O)]^2}{E_{\lambda,d}\zeta_t^2(O)}\geq\frac{C^2(\lambda,d,\zeta)t^{-d}}{F_t(O)}.
\end{align}

We denote by $\|\cdot\|_{\infty}$ the $l_{\infty}$ norm on
$R^{\mathbb{Z}^d}$ such that
\[
\|\zeta\|_{\infty}=\sup_{x\in \mathbb{Z}^d}|\zeta(x)|
\]
for any $\zeta\in R^{\mathbb{Z}^d}$.

By direct calculation, for any $\zeta_1,\zeta_2\in R^{\mathbb{Z}^d}$
such that $\|\zeta_1\|_{\infty},\|\zeta_2\|_{\infty}<+\infty$,
\begin{equation}\label{equ 4.7}
\|G\zeta_1-G\zeta_2\|_{\infty}\leq (1+8\lambda d+4\lambda
d^2)\|\zeta_1-\zeta_2\|_{\infty}.
\end{equation}
By $\eqref{equ 4.7}$ and classic Theory for ODE on Banach Space, ODE
\eqref{equ 4.2 two} has the unique solution such that
\begin{equation}\label{equ 4.8}
F_t=\Gamma_tF_0,
\end{equation}
where
\begin{equation}\label{equ 4.9}
\Gamma_t=\big(\gamma_t(x,y)\big)_{x\in \mathbb{Z}^d,y\in
\mathbb{Z}^d}=\exp\{tG\}=\sum_{n=0}^{+\infty}\frac{t^nG^n}{n!}.
\end{equation}
\eqref{equ 4.7} ensures the sum in \eqref{equ 4.9} is finite. By
\eqref{equ 4.8},
\begin{equation}\label{equ 4.10}
F_t(O)=\sum_{x\in \mathbb{Z}^d}\gamma_t(O,x)H(x),
\end{equation}
since $F_0=H$.

Since $H$ is the eigenvector of $G$ with respect to the eigenvalue
$\mu$, $H$ is also the eigenvector of $\Gamma_t=e^{tG}$ with respect
to the eigenvalue $\exp\{t\mu\}$.

As a result,
\begin{equation}\label{equ 4.12}
F_t(O)=\sum_{x\in \mathbb{Z}^d}\gamma_t(O,x)H(x)=e^{t\mu}H(O).
\end{equation}

By \eqref{equ 4.6} and \eqref{equ 4.12},
\[
P_{\lambda,d}(\eta_t^{\delta_1}(O)=1)\geq
\frac{C^2(\lambda,d,\zeta)t^{-d}}{H(O)}e^{-t\mu},
\]
and
\[
I(\lambda,d)=\lim_{t\rightarrow+\infty}\frac{1}{t}\log
P_{\lambda,d}(\eta_t^{\delta_1}(O)=1)\geq -\mu.
\]

\qed

To search $\lambda$ and $\mu$ satisfies the condition in Lemma
\ref{lemma 4.4}, we introduce the simple random walk on
$\mathbb{Z}^d\cup \{\triangle\}$, where $\triangle\not\in Z^d$ is an
absorbed state.

For $d\geq 1$ and $p\in [0,1]$, let $\{S_n(d,p):n=0,1,2,\ldots\}$ be
simple random walk on $\mathbb{Z}^d\cup \{\triangle\}$ with
transition probability
\begin{align}\label{equ 4.13 random walk tran function}
&P\big(S_{n+1}(d,p)=y\big|S_n(d,p)=x\big)=\frac{p}{2d},\notag\\
&P\big(S_{n+1}(d,p)=\triangle\big|S_n(d,p)=x\big)=1-p,\\
&P\big(S_{n+1}(d,p)=\triangle\big|S_n(d,p)=\triangle\big)=1 \notag
\end{align}
for $n\geq 0$, each $x\in Z^d$ and each $y\sim x$.

For $d\geq 1$, $p\in [0,1]$ and $x\in \mathbb{Z}^d$, we define
\[
\tau(d,p)=\inf\{n\geq 0:S_n(d,p)=O\}
\]
and
\[
R(x,d,p)=P\big(\tau(d,p)<+\infty\big|S_0(d,p)=x\big).
\]
We will give $H(x)$ with the form $R(x,d,p)$. For this purpose, we
need the following lemma.
\begin{lemma}\label{lemma 4.5} For $d\geq
1$ and each $x\in \mathbb{Z}^d$, $R(x,d,p)$ is continuous in $p$.
\end{lemma}

\proof

The conclusion is trivial for $x=O$. For $x\neq O$ and $0\leq
p_1<p_2\leq 1$, we construct a coupling for $\{S_n(d,p_1)\}_{n\geq
0}$ and $\{S_n(d,p_2)\}_{n\geq 0}$ with $S_0(d,p_1)=S_0(d,p_2)=x$
such that $S_n(d,p_1)=\triangle$ or $S_n(d,p_1)=S_n(d,p_2)\neq
\triangle$ for each $n\geq 1$.

The transition probability matrix $\widehat{P}$ of the coupling
process $\{S_n(d,p_1),S_n(d,p_2)\}$ is given by
\begin{equation}\label{equ 4.14}
\widehat{P}\big((x_1,y_1),(x_2,y_2)\big)=
\begin{cases}
\frac{p_1}{2d} &\text{~if~}x_1=y_1\neq \triangle, x_2=y_2\sim x_1,\\
1-p_2 &\text{~if~}x_1=y_1\neq \triangle, x_2=y_2=\triangle,\\
\frac{p_2-p_1}{2d} &\text{~if~}x_1=y_1\neq \triangle, x_2=\triangle,y_2\sim y_1,\\
\frac{p_2}{2d} &\text{~if~}x_1=\triangle, y_1\neq \triangle,
x_2=\triangle,y_2\sim y_1,\\
1-p_2 &\text{~if~}x_1=\triangle, y_1\neq \triangle,
x_2=y_2=\triangle,\\
0&\text{~otherwise}.
\end{cases}
\end{equation}
It is easy to check that $\widehat{P}$ gives a coupling of
$S_n(d,p_1)$ and $S_n(d,p_2)$ by direct calculation. The coupling
ensures that $S_n(d,p_2)=S_n(d,p_1)\neq \triangle$ when
$S_n(d,p_1)\neq \triangle$.

As a result, conditioned on $S_0(d,p_1)=S_0(d,p_2)=x$,
\begin{align}\label{equ 4.15}
R(x,d,p_2)-R(x,d,p_1)&=P\big(\tau(d,p_1)=+\infty,\tau(d,p_2)<+\infty\big)\notag\\
&\leq P\big(\exists~n>0,S_n(d,p_1)=\triangle,S_n(d,p_2)\neq
\triangle\big).
\end{align}
Let
\[
\beta=\inf\{n\geq 1: S_n(d,p_1)=\triangle,S_n(d,p_2)\neq
\triangle\},
\]
then, by \eqref{equ 4.14} and \eqref{equ 4.15},
\begin{align}\label{equ 4.16}
&R(x,d,p_2)-R(x,d,p_1)\leq P(\beta<+\infty)
=\sum_{l=1}^{+\infty}P(\beta=l)\notag\\
&=\sum_{l=1}^{+\infty}P(\beta>l-1,S_l(d,p_1)=\triangle,S_l(d,p_2)\neq
\triangle)\notag\\
&=\sum_{l=1}^{\infty}P(S_{l-1}(d,p_1)\neq \triangle)(p_2-p_1)\notag\\
&=\sum_{l=1}^{\infty}p_1^{l-1}(p_2-p_1)=\frac{p_2-p_1}{1-p_1}.
\end{align}
Lemma \ref{lemma 4.5} follows from \eqref{equ 4.16}.

\qed

Now we give a lower bound for $I(\lambda,d)$.
\begin{lemma}\label{lemma 4.6}
For each $d\geq 1$ and $\lambda<\frac{1}{2d}$, there exists unique
$p=p(\lambda,d)\in (0,1)$ such that
\begin{equation}\label{equ 4.17}
1+2\lambda d R(e_1,d,p)=\frac{4\lambda d}{p}\big[1-dR(e,d,p)\big].
\end{equation}
Furthermore,
\begin{equation}\label{equ 4.18}
I(\lambda,d)\geq -4\lambda d\big[\frac{1}{p(\lambda,d)}-1\big]
\end{equation}
for $\lambda<\frac{1}{2d}$.
\end{lemma}

\proof

For $p\in (0,1]$, we define
\[
K(p)=\frac{4\lambda d}{p}\big[1-dR(e_1,d,p)\big]-1-2\lambda
dR(e_1,d,p).
\]
It is obviously that $K(p)$ is decreasing in $p$. By Lemma
\ref{lemma 4.5}, $K(p)$ is continuous in $p$.

Since $\lambda<\frac{1}{2d}$ and  $R(e_1,d,1)\geq
P\big(S_1(d,1)=O\big|S_0(d,1)=e_1\big)=1/2d$,
\begin{equation}\label{equ 4.19}
K(1)<0.
\end{equation}
Since $R(e_1,d,0)=0$,
\begin{equation}\label{equ 4.20}
\lim_{p\rightarrow 0}K(p)=+\infty.
\end{equation}
The existence and uniqueness of $p(\lambda,d)$ follows from
\eqref{equ 4.19}, \eqref{equ 4.20} and the fact that $K(p)$ is
continuous and decreasing in $p$.

\quad

Let $\mu=4\lambda d\big[1/p(\lambda,d)-1\big]$,
$H(x)=R\big(x,d,p(\lambda,d)\big)$ for each $x\in \mathbb{Z}^d$,
then according to the fact that $p(\lambda,d)$ satisfies \eqref{equ
4.17} and
\[
R(x,d,p)=\frac{p}{2d}\sum_{y:y\sim x}R(y,d,p)
\]
for each $x\neq O$, it is easy to check that
\[
GH=\mu H.
\]
As a result, \eqref{equ 4.18} follows from Lemma \ref{lemma 4.4}.

\qed

To give a limit theorem of $p(\lambda,d)$, we need the following
lemma.
\begin{lemma}\label{lemma 4.7} For $\{p_d\}_{d=1,2,\ldots}$ such
that $p_d\in (0,1)$ for each $d\geq 1$, if
\[
\lim_{d\rightarrow+\infty}p_d=c,
\]
then
\[
\lim_{d\rightarrow+\infty}2dR(e_1,d,p_d)=c.
\]
\end{lemma}

\proof
\begin{align}\label{equ 4.21}
R(e_1,d,p_d)&\geq P\big(\tau(d,p_d)=1\big|S_0(d,p_d)=e_1\big)\notag\\
&=P\big(S_1(d,p_d)=O\big|S_0(d,p_d)=e_1\big)=\frac{p_d}{2d}.
\end{align}
On the other hand,
\begin{align}\label{equ 4.22}
R(e_1,d,p_d)&=\frac{p_d}{2d}+P\big(2\leq\tau(d,p_d)<+\infty\big|S_0(d,p_d)=e_1\big)\notag\\
&\leq
\frac{p_d}{2d}+P\big(2\leq\tau(d,1)<+\infty\big|S_0(d,1)=e_1\big).
\end{align}
According to Lemma 5.3 of \cite{Xue2014},
\begin{equation}\label{equ 4.23}
\lim_{d\rightarrow+\infty}dP\big(2\leq\tau(d,1)<+\infty\big|S_0(d,1)=e_1\big)=0.
\end{equation}
Lemma \ref{lemma 4.7} follows from \eqref{equ 4.21},\eqref{equ 4.22}
and \eqref{equ 4.23}.

\qed

\quad

Finally, we give the proof of
$\liminf_{d\rightarrow+\infty}I(\lambda/d,d)\geq 2\lambda-1$ for
$\lambda\in (0,1/2)$.

\proof

For any $\lambda\in (0,1/2)$, we define
\[
\overline{c}(\lambda)=\limsup_{d\rightarrow+\infty}p(\lambda/d,d)
\]
and
\[
\underline{c}(\lambda)=\liminf_{d\rightarrow+\infty}p(\lambda/d,d).
\]
Then by \eqref{equ 4.17} and Lemma \ref{lemma 4.7},
\[
1=\frac{4\lambda}{\overline{c}(\lambda)}\big[1-\frac{\overline{c}(\lambda)}{2}\big]
\]
and
\[
1=\frac{4\lambda}{\underline{c}(\lambda)}\big[1-\frac{\underline{c}(\lambda)}{2}\big].
\]
Therefore,
\[
\overline{c}(\lambda)=\underline{c}(\lambda)=c(\lambda)=\frac{4\lambda}{1+2\lambda}
\]
and hence
\begin{equation}\label{equ 4.24}
\lim_{d\rightarrow+\infty}p(\lambda/d,d)=c(\lambda).
\end{equation}
By \eqref{equ 4.18} and \eqref{equ 4.24},
\[
\liminf_{d\rightarrow+\infty}I(\lambda/d,d)\geq
-4\lambda[\frac{1}{c(\lambda)}-1]=2\lambda-1
\]
for $\lambda\in (0,1/2)$.

\qed

To finish the proof of Theorem \ref{theorem 2.2 main}, we only need
to deal with the case where $\lambda=1/2$.

\proof[Proof of Theorem \ref{theorem 2.2 main}]

For $\lambda\in (0,1/2)$, we have shown that
\[
\limsup_{d\rightarrow+\infty}I(\lambda/d,d)\leq 2\lambda-1 \] in
Section \ref{section 3 upper bound} and
\[
\liminf_{d\rightarrow+\infty}I(\lambda/d,d)\geq 2\lambda-1
\]
in this section. Therefore,
\begin{equation}\label{equ 4.25}
\lim_{d\rightarrow+\infty}I(\lambda/d,d)=2\lambda-1
\end{equation}
for $\lambda\in (0,1/2)$.

In Section \ref{section 3 upper bound}, we show that
\begin{equation}\label{equ 4.26}
\lim_{d\rightarrow+\infty}I(\lambda/d,d)=0
\end{equation} for
$\lambda>1/2$. It is obviously that $I(\lambda,d)$ is increasing in
$\lambda$. Therefore, by \eqref{equ 4.25} and \eqref{equ 4.26},
\begin{equation}\label{equ 4.27}
\lim_{d\rightarrow+\infty}I(1/2d,d)=0.
\end{equation}
Theorem \ref{theorem 2.2 main} follows from \eqref{equ 4.25},
\eqref{equ 4.26} and \eqref{equ 4.27}.

\qed

Now the whole proof of Theorem \ref{theorem 2.2 main} is completed.
Furthermore, we show that
\[
-4\lambda d\big[\frac{1}{p(\lambda,d)}-1\big]\leq I(\lambda,d)\leq
2\lambda d-1
\]
for $\lambda<1/2d$, where $p(\lambda,d)$ is the unique solution to
\[
1+2\lambda d R(e_1,d,p)=\frac{4\lambda d}{p}\big[1-dR(e,d,p)\big].
\]

\quad

We give a sketch for the proof of Theorem \ref{theorem 2.3}.

\proof[Proof of \ref{theorem 2.3}]

We only need to show that
$\liminf_{d\rightarrow+\infty}J(\lambda/d,d)\geq 2\lambda-1$ for
$\lambda\in (0,1/2)$.

Let $\{\alpha_t\}_{t\geq 0}$ be Markov processes with state space
$[0,+\infty)^{\mathbb{Z}^d}$ such that the generator
$\widetilde{\Omega}$ of $\{\alpha_t\}_{t\geq 0}$ is given by
\begin{align*}
\widetilde{\Omega}f(\alpha)&=\sum_{x\in \mathbb{Z}^d}\big[f(\alpha^{x,0})-f(\alpha)\big]\\
&+\lambda\sum_{x\in \mathbb{Z}^d}\sum_{y:y\sim
x}\big[f(\alpha^{x,\alpha(x)+\alpha(y)})-f(\alpha)\big] +\sum_{x\in
\mathbb{Z}^d}(1-2\lambda d)f^{\prime}_x(\alpha)\alpha(x),
\end{align*}
where
\[
\alpha^{x,m}(y)=
\begin{cases}
\alpha(y) &\text{~if~} y\neq x,\\
m &\text{~if~} y=x
\end{cases}
\]
for $x\in \mathbb{Z}^d$ and $m\geq 0$.

When $\alpha_0(x)>0$ for each $x\in \mathbb{Z}^d$, then according to
a similar analysis with that in the proof of Lemma \ref{lemma 4.1},
\[
\beta_t^{\delta_1}(O)=
\begin{cases}
1&\text{~if~} \alpha_t(O)=1,\\
0&\text{~if~} \alpha_t(O)=0
\end{cases}
\]
in the sense of coupling and hence
\begin{equation}\label{equ 4.28}
P_{\lambda,d}(\beta^{\delta_1}_t(O)=1)=P_{\lambda,d}(\alpha_t(O)>0)\geq
\frac{\big[E_{\lambda,d}\alpha_t(O)\big]^2}{E_{\lambda,d}\alpha_t^2(O)}.
\end{equation}
We define $\widetilde{F}:[0,+\infty)\rightarrow
[0,+\infty)^{\mathbb{Z}^d}$ such that
\[
\widetilde{F}_t(x)=E_{\lambda,d}\big[\alpha_t(O)\alpha_t(x)\big]
\]
for $x\in \mathbb{Z}^d$ and $t\geq 0$. Then,
\[
\frac{d}{dt}\widetilde{F}_t=\widetilde{G}\widetilde{F}_t,
\]
where $\widetilde{G}$ is a $\mathbb{Z}^d\times \mathbb{Z}^d$ matrix
such that
\[
\widetilde{G}(x,y)=
\begin{cases}
2\lambda &\text{~if~} x\neq O, y\sim x,\\
-4\lambda d &\text{~if~} x\neq O, y=x,\\
1-2\lambda d &\text{~if~} x=y=O,\\
4\lambda d &\text{~if~} x=O, y=e_1,\\
0 &\text{~otherwise}.
\end{cases}
\]
When $\lambda<\frac{1}{2d}$, according to a similar analysis with
that in the proof of Lemma \ref{lemma 4.6}, there exists unique
$\widetilde{p}=\widetilde{p}(\lambda,d)$ such that
\[
\frac{4\lambda d}{\widetilde{p}}-2\lambda d=1+4\lambda d
R(e_1,d,\widetilde{p}).
\]
Let $\widetilde{\mu}=4\lambda d\big[1/\widetilde{p}-1\big]$ and
$\widetilde{H}(x)=R(x,d,\widetilde{p})$ for each $x\in
\mathbb{Z}^d$, then
\[
\widetilde{G}\widetilde{H}=\widetilde{\mu}\widetilde{H}.
\]
According to \eqref{equ 4.28} and a similar analysis with that in
the proof of Lemma \ref{lemma 4.4},
\[
J(\lambda,d)\geq -\widetilde{\mu}
\]
for $\lambda<\frac{1}{2d}$.

Since $R(e_1,d,\widetilde{p})\leq R(e_1,d,1)\rightarrow 0$ as
$d\rightarrow+\infty$,
\[
\lim_{d\rightarrow+\infty}\widetilde{p}(\lambda/d,d)=\frac{4\lambda}{1+2\lambda}
\]
for $\lambda<1/2$. As a result,
\[
\liminf_{d\rightarrow+\infty}J(\lambda/d,d)\geq
-4\lambda\big[\frac{1+2\lambda}{4\lambda}-1\big]=2\lambda-1
\]
for $\lambda<1/2$.

\qed

\quad

\textbf{Acknowledgments.} The author is grateful to the financial
support from the National Natural Science Foundation of China with
grant number 11171342.

{}

\begin{thebibliography}{}
\bibitem{And1992}Andjel, E. D., Liggett, T. M. and Mountford, T. (1992). Clustering in one-dimensional threshold voter models. \emph{Stochastic Processes and Their Applications} \textbf{42}, 73-90.
\bibitem{Dur1991}Cox, J. T. and Durrett, R. (1991). Nonlinear voter models. In \emph{Random Walks, Brownian Motion
and Interacting Particle Systems. A Festschrift in Honor of Frank Spiter} 189-201. Birkh\"{a}user, Boston.
\bibitem{Fonte2008}Fontes, L. R., Schonmann, R. H. (2008). Threshold \(\theta\geq2\) contact processes on homogeneous
trees. \emph{Probability Theory and Related Fields} \textbf{141}, 513-541.
\bibitem{Grif1983}Griffeath, D. (1983). The Binary Contact Path Process. \emph{The Annals of Probability} \textbf{11} 692-705.
\bibitem{Handjani1999}Handjani, S. (1999). The complete convergence theorem for coexistent threshold voter models. \emph{The Annals of Probability} \textbf{27} 226-245.
\bibitem{LIG1985}Liggett, T. M. (1985). \emph{Interacting Particle Systems.} Springer, New York.
\bibitem{LIG1994}Liggett, T. M. (1994). Coexistence in threshold voter models. \emph{The Annals of Probability}. \textbf{22}, 764-802.
\bibitem{LIG1999}Liggett, T. M. (1999). \emph{Stochastic interacting systems: contact, voter and exclusion processes.}
Springer, New York.
\bibitem{Mou2009}Mountford, T. and Schonmann, R. H. (2009) The survival of large dimensional threshold contact processes. \emph{The Annals of Probability} \textbf{37}, 1483-1501.
\bibitem{Xue2012}Xue, XF. (2012). Critical density points for threshold voter models on homogeneous trees. \emph{Journal of Statistical Physics}. \textbf{146}, 423-433.
\bibitem{Xue2014}Xue, XF. (2014). Asymptotic behavior of critical infection rates for threshold-one contact processes on lattices and regular trees. \emph{Journal of Theoretical Probability}. Published online on February 2014.
\bibitem{Xue2015}Xue, XF. (2015). Fluid limit of threshold voter models on tori. \emph{Journal of Statistical Physics}. Published online on January  2015.
\end{thebibliography}
\end{document}